\documentclass[12pt,thmsa]{article}
\usepackage[latin9]{inputenc}
\usepackage[a4paper]{geometry}
\usepackage{varioref}
\usepackage{url}
\usepackage{amsmath}
\usepackage{amssymb}
\usepackage{setspace}
\usepackage[authoryear]{natbib}
\PassOptionsToPackage{normalem}{ulem}
\usepackage{ulem}
\usepackage[unicode=true,
 bookmarks=false,
 breaklinks=false,pdfborder={0 0 1},backref=section,colorlinks=false]
 {hyperref}

\makeatletter

\providecommand{\tabularnewline}{\\}

\usepackage{amsfonts}
\usepackage{url}
\usepackage{amsthm}
\usepackage{pdflscape}
\usepackage{graphicx}
\usepackage{url}
\graphicspath{{/Users/bgraham/Dropbox/Research/Networks/DyadicRegression/Graphics/}}

\makeatother

\begin{document}
\begin{singlespacing}
\title{Kernel Density Estimation for Undirected Dyadic Data}

\maketitle
\medskip{}

\noindent \begin{center}
{\large{}{}Bryan S. Graham}\footnote{{\footnotesize{}{}Department of Economics, University of California
- Berkeley, 530 Evans Hall \#3380, Berkeley, CA 94720-3880 and National
Bureau of Economic Research, {}}{\footnotesize{}\uline{e-mail:}}{\footnotesize{}{}
\href{http://bgraham\%5C\%40econ.berkeley.edu}{bgraham@econ.berkeley.edu},
{}}{\footnotesize{}\uline{web:}}{\footnotesize{}{} \url{http://bryangraham.github.io/econometrics/}.
}\\
 {\large{}{}$^{\#}$}{\footnotesize{}{}Department of Economics,
University of California - Berkeley, 530 Evans Hall \#3380, Berkeley,
CA 94720-3880, {}}{\footnotesize{}\uline{e-mail:}}{\footnotesize{}{}
\href{http://fniu\%5C\%40berkeley.edu}{fniu@berkeley.edu}.}{\large{}{}}\\
 {\large{}{}$^{\dagger}$}{\footnotesize{}{}Department of Economics,
University of California - Berkeley, 530 Evans Hall \#3380, Berkeley,
CA 94720-3880, {}}{\footnotesize{}\uline{e-mail:}}{\footnotesize{}{}
\href{http://powell\%5C\%40econ.berkeley.edu}{powell@econ.berkeley.edu}.\medskip{}
 }\\
 {\footnotesize{}{}We thank Michael Jansson and Konrad Menzel for
helpful conversation as well as audiences at Berkeley, Brown, Toulouse,
Warwick, Bristol, University College London and the Conference Celebrating
Whitney Newey's Contributions to Econometrics for useful questions
and suggestions. All the usual disclaimers apply. Financial support
from NSF grant SES \#1851647 is gratefully acknowledged.}}{\large{}{}, Fengshi Niu$^{\#}$ and James L. Powell$^{\dagger}$}{\large\par}
\par\end{center}

\begin{center}
 
\par\end{center}

\noindent \begin{center}
\medskip{}
 \textsc{\large{}{}}\textsc{\large{}\uline{Initial Draft}}\textsc{\large{}{}:
January 2019, {}}\textsc{\large{}\uline{This Draft}}\textsc{\large{}{}:
July 2019}{\large\par}
\par\end{center}

\begin{center}
 
\par\end{center}

\pagebreak{} 
\begin{center}
{\Large{}{}Abstract}{\Large\par}
\par\end{center}

\begin{center}
 
\par\end{center}

We study nonparametric estimation of density functions for undirected
dyadic random variables (i.e., random variables defined for all $n\overset{def}{\equiv}\tbinom{N}{2}$
unordered pairs of agents/nodes in a weighted network of order $N$).
These random variables satisfy a local dependence property: any random
variables in the network that share one or two indices may be dependent,
while those sharing no indices in common are independent. In this
setting, we show that density functions may be estimated by an application
of the kernel estimation method of \citet{Rosenblatt_AMS56} and \citet{Parzen_AMS62}.
We suggest an estimate of their asymptotic variances inspired by a
combination of (i) Newey's \citeyearpar{Newey_ET94b} method of variance
estimation for kernel estimators in the ``monadic'' setting and
(ii) a variance estimator for the (estimated) density of a simple
network first suggested by \citet{Holland_Leinhardt_SM76}. More unusual
are the rates of convergence and asymptotic (normal) distributions
of our dyadic density estimates. Specifically, we show that they converge
at the same rate as the (unconditional) dyadic sample mean: the square
root of the number, $N$, of nodes. This differs from the results
for nonparametric estimation of densities and regression functions
for monadic data, which generally have a slower rate of convergence
than their corresponding sample mean.

{\footnotesize{}{}\medskip{}
}{\footnotesize\par}

\qquad{}\textbf{JEL Classification:} C24, C14, C13.

{\footnotesize{}{}\medskip{}
 }{\footnotesize\par}

\qquad{}\textbf{Keywords:} Networks, Dyads, Kernel Density Estimation

\end{singlespacing}

\thispagestyle{empty} \pagebreak{}

\setcounter{page}{1}

\section{Introduction}

Many important social and economic variables are naturally defined
for pairs of agents (or dyads). Examples include trade between pairs
of countries \citep[e.g.,][]{Tinbergen_SWE62}, input purchases and
sales between pairs of firms \citep[e.g.,][]{Atalay_et_al_PNAS11},
research and development (R\&D) partnerships across firms \citep[e.g.,][]{Konig_et_al_RESTAT18}
and friendships between individuals \citep[e.g.,][]{Christakis_et_al_NBER10}.
Dyadic data arises frequently in the analysis of social and economic
networks. In economics such analyses are predominant in, for example,
the analysis of international trade flows. See \citet{Graham_HBE18}
for many other examples and references.

While the statistical analysis of network data began almost a century
ago, rigorously justified methods of inference for network statistics
are only now emerging \citep[cf.,][]{Goldenberg_etal_FTML09}. In
this paper we study nonparametric estimation of the density function
of a (continuously-valued) dyadic random variable. Examples included
the density of migration across states, trade across nations, liabilities
across banks, or minutes of telephone conversation among individuals.
While nonparametric density estimation using independent and identically
distributed random samples, henceforth ``monadic'' data, is well-understood,
its dyadic counterpart has, to our knowledge, not yet been studied.

\citet{Holland_Leinhardt_SM76} derived the sampling variance of the
link frequency in a simple network (and of other low order subgraph
counts). A general asymptotic distribution theory for subgraph counts,
exploiting recent ideas from the probability literature on dense graph
limits \citep[e.g.,][]{Diaconis_Janson_RM08,Lovasz_AMS12}, was presented
in \citet{Bickel_et_al_AS11}.\footnote{See \citet{Nowicki_SN91} for a summary of earlier research in this
area.} \citet{Menzel_arXiv17} presents bootstrap procedures for inference
on the mean of a dyadic random variable. Our focus on nonparametric
density estimation appears to be novel. Density estimation is, of
course, a topic of intrinsic interest to econometricians and statisticians,
but it also provides a relatively simple and canonical starting point
for understanding nonparametric estimation more generally. In the
conclusion of this paper we discuss ongoing work on other non- and
semi-parametric estimation problems using dyadic data.

We show that an (obvious) adaptation of the \citet{Rosenblatt_AMS56}
and \citet{Parzen_AMS62} kernel density estimator is applicable to
dyadic data. While our dyadic density estimator is straightforward
to define, its rate-of-convergence and asymptotic sampling properties,
depart significantly from its monadic counterpart. Let $N$ be the
number of sampled agents and $n=\tbinom{N}{2}$ the corresponding
number of dyads. Estimation is based upon the $n$ dyadic outcomes.
Due to dependence across dyads sharing an agent in common, the rate
of convergence of our density estimate is (generally) much \emph{slower}
than it would be with $n$ i.i.d. outcomes. This rate-of-convergence
is also invariant across a wide range of bandwidth sequences. This
property is familiar from the econometric literature on semiparametric
estimation \citep[e.g.,][]{Powell_HBE94}. Indeed, from a certain
perspective, our nonparametric dyadic density estimate can be viewed
as a semiparametric estimator (in the sense that it can be thought
of as an average of nonparametrically estimated densities). We also
explore the impact of ``degeneracy'' -- which arises when dependence
across dyads vanishes -- on our sampling theory; such degeneracy
features prominently in Menzel's \citeyearpar{Menzel_arXiv17} innovative
analysis of inference on dyadic means. We expect that many of our
findings generalize to other non- and semi-parametric network estimation
problems.

In the next section we present our maintained data/network generating
process and proposed kernel density estimator. Section \ref{sec: MSE}
explores the mean square error properties of this estimator, while
Section \ref{sec: asymptotic} outlines asymptotic distribution theory.
Section \ref{sec: asymptotic_variance_estimation} presents a consistent
variance estimator, which can be used to construct Wald statistics
and Wald-based confidence intervals. We summarize the results of a
small simulation study in Section \ref{sec: simulation_study}. In
Section \ref{sec: extensions} we discuss various extensions and ongoing
work. Calculations not presented in the main text are collected in
Appendix \ref{app: Proofs}.

It what follows we interchangeably use unit, node, vertex, agent and
individual all to refer to the $i=1,\ldots,N$ vertices of the sampled
network or graph. We denote random variables by capital Roman letters,
specific realizations by lower case Roman letters and their support
by blackboard bold Roman letters. That is $Y$, $y$ and $\mathbb{Y}$
respectively denote a generic random draw of, a specific value of,
and the support of, $Y$. For $W_{ij}$ a dyadic outcome, or weighted
edge, associated with agents $i$ and $j$, we use the notation $\mathbf{W}=\left[W_{ij}\right]$
to denote the $N\times N$ adjacency matrix of all such outcomes/edges.
Additional notation is defined in the sections which follow.

\section{\label{sec: model}Model and estimator}

\subsubsection*{Model}

Let $i=1,\ldots,N$ index a simple random sample of $N$ agents from
some large (infinite) network of interest. A pair of agents constitutes
a \emph{dyad}. For each of the $n=\tbinom{N}{2}$ sampled dyads, that
is for $i=1,...,N-1$ and $j=i+1,\ldots,N$, we observe the (scalar)
random variable $W_{ij}$, generated according to 
\begin{equation}
W_{ij}=W(A_{i},A_{j},V_{ij})=W(A_{j},A_{i},V_{ij}),\label{eq: DGP}
\end{equation}
where $A_{i}$ is a node-specific random vector of attributes (of
arbitrary dimension, not necessarily observable), and $V_{ij}=V_{ji}$
is an unobservable scalar random variable which is continuously distributed
on $\mathbb{R}$ with density function $f_{V}(v)$.\footnote{In words we observe the weighted subgraph induced by the randomly
sampled agents.} Observe that the function $W(a_{1},a_{2},v_{12})$ is symmetric in
its first two arguments, ensuring that $W_{ij}=W_{ji}$ is undirected.

In what follows we directly maintain \eqref{eq: DGP}, however, it
also a consequence of assuming that the infinite graph sampled from
is jointly exchangeable \citep{Aldous_JMA81,Hoover_WP79}. Joint exchangeability
of the sampled graph $\mathbf{W}=\left[W_{ij}\right]$ implies that
\begin{equation}
\left[W_{ij}\right]\overset{D}{=}\left[W_{\pi\left(i\right)\pi\left(j\right)}\right]\label{eq: exchangeability}
\end{equation}
for every $\pi\in\Pi$ where $\pi:\left\{ 1,\ldots,N\right\} \rightarrow\left\{ 1,\ldots,N\right\} $
is a permutation of the node indices. Put differently, when node labels
have no meaning we have that the ``likelihood'' of any simultaneous
row and column permutation of $\mathbf{W}$ is the same as that of
$\mathbf{W}$ itself.\footnote{For $\mathbf{W}=\left[W_{ij}\right]$ the $N\times N$ weighted adjacency
matrix and $\mathbf{P}$ any conformable permutation matrix 
\[
\Pr\left(\mathbf{W}\leq\mathbf{w}\right)=\Pr\left(\mathbf{P}\mathbf{W}\mathbf{P}\leq\mathbf{w}\right)
\]
for all $\mathbf{w}\in\mathbb{W=R}^{\tbinom{N}{2}}.$} See \citet{Menzel_arXiv17} for a related discussion.

Our target object of estimation is the marginal density function $f_{W}(w)$
of $W_{ij}$, defined as the derivative of the cumulative distribution
function (c.d.f.) of $W_{ij},$ 
\[
\Pr\{W_{ij}\leq w\}\overset{def}{\equiv}F_{W}(w)=\int_{-\infty}^{w}f_{W}(u)\mathrm{d}u.
\]
To ensure this density function is well-defined on the support of
$W_{ij},$ we assume that the unknown function $W(a_{1},a_{2},v)$
is strictly increasing and continuously differentiable in its third
argument $v$, and we also assume that $A_{i}$ and $A_{j}$ are statistically
independent of the ``error term'' $V_{ij}$ for all $i$ and $j.$
Under these assumptions, by the usual change-of-variables formula,
the conditional density of $W_{ij}$ given $A_{i}=a_{1}$ and $A_{j}=a_{2}$
takes the form 
\[
f_{Y|AA}(w|a_{1},a_{2})=f_{V}(W^{-1}(a_{1},a_{2},w))\cdot\left\vert \frac{\partial W(a_{1},a_{2},W^{-1}(a_{1},a_{2},w))}{\partial v}\right\vert ^{-1}.
\]
In the derivations below we will assume this density function is bounded
and twice continuously differentiable at $w$ with bounded second
derivative for all $a_{1}$ and $a_{2}$; this will follow from the
similar smoothness conditions imposed on the primitives $W^{-1}(\cdot,\cdot,w)$
and $f_{V}(v).$

To derive the marginal density of $W_{ij}$ note that, by random sampling,
the $\{A_{i}\}$ sequence is independently and identically distributed
(i.i.d.), as is the $\{V_{ij}\}$ sequence. Under these conditions,
we can define the conditional densities of $W_{ij}$ given $A_{i}=a$
or $A_{j}=a$ alone as 
\[
f_{W|A}(w|a)\equiv\mathbb{E}[f_{W|AA}(w|a,A_{j})]=\mathbb{E}[f_{W|AA}(w|A_{i},a)],
\]
and, averaging, the marginal density of interest as 
\[
f_{W}(w)\overset{def}{\equiv}\mathbb{E}[f_{W|AA}(w|A_{i},A_{j})]=\mathbb{E}[f_{W|A}(w|A_{i})].
\]

Let $i,j,k$ and $l$ index distinct agents. The assumption that $\{A_{i}\}$
and $\{V_{ij}\}$ are i.i.d. implies that while $W_{ij}$ varies independently
of $W_{kl}$ (since the $\left\{ i,j\right\} $ and $\left\{ k,l\right\} $
dyads share no agents in common), $W_{ij}$ will not vary independently
of $W_{ik}$ as both vary with $A_{i}$ (since the $\left\{ i,j\right\} $
and $\left\{ i,k\right\} $ dyads both include agent $i$). This type
of dependence structure is sometimes referred to as ``dyadic clustering''
in empirical social science research \citep[cf.,][]{Fafchamp_Gubert_JDE07,Cameron_Miller_WP14,Aronow_et_al_PA17}.
The implications of this dependence structure for density estimation
and -- especially -- inference is a key area of focus in what follows.

\subsubsection*{Estimator}

Given this construction of the marginal density $f_{W}(w)$ of $W_{ij},$
it can be estimated using an immediate extension of the kernel density
estimator for monadic data first proposed by \citet{Rosenblatt_AMS56}
and \citet{Parzen_AMS62}: 
\begin{align*}
\hat{f}_{W}(w) & =\tbinom{N}{2}^{-1}\sum_{i=1}^{N-1}\sum_{j=1+1}^{N}\frac{1}{h}K\left(\frac{w-W_{ij}}{h}\right)\\
 & \overset{def}{\equiv}\frac{1}{n}\sum_{i<j}K_{ij},
\end{align*}
where 
\[
K_{ij}\overset{def}{\equiv}\frac{1}{h}K\left(\frac{w-W_{ij}}{h}\right).
\]
Here $K(\cdot)$ is a kernel function assumed to be (i) bounded ($K(u)\leq\bar{K}$
for all $u$), (ii) symmetric ($K(u)=K(-u)$), (ii) , and zero outside
a bounded interval ($K(u)=0$ if $\left\vert u\right\vert >\bar{u}$);
we also require that it (iv) integrates to one ($\int K(u)du=1$).
The bandwidth $h=h_{N}$ is assumed to be a positive, deterministic
sequence (indexed by the number of nodes $N$) that tends to zero
as $N\rightarrow\infty,$ and will satisfy other conditions imposed
below. A discussion of the motivation for the kernel estimator $\hat{f}_{W}(w)$
and its statistical properties under random sampling (of monadic variables)
can be found in \citet[Chapers 2 \& 3]{Silverman_DESDA86}.

\section{\label{sec: MSE}Rate of convergence analysis}

To formulate conditions for consistency of $\hat{f}_{W}(w),$ we will
evaluate its expectation and variance, which will yield conditions
on the bandwidth sequence $h_{N}$ for its mean squared error to converge
to zero.

A standard calculation yields a bias of $\hat{f}_{W}(w)$ equal to
(see Appendix \ref{app: Proofs}) 
\begin{align}
E\left[\hat{f}_{W}(w)\right]-f_{W}(w) & =h^{2}B(w)+o(h^{2})\label{eq: bias}\\
 & =O(h_{N}^{2}),\nonumber 
\end{align}
with 
\[
B\left(w\right)\overset{def}{\equiv}\frac{1}{2}\frac{\partial^{2}f_{W}(w)}{\partial w{}^{2}}\int u^{2}K\left(u\right)\mathrm{d}u.
\]
Equation \eqref{eq: bias} coincides with the bias of the kernel density
estimate based upon a random (``monadic'') sample.

The expression for the variance of $\hat{f}_{W}(w)$, in contrast
to that for bias, does differ from the monadic (i.i.d.) case due to
the (possibly) nonzero covariance between $K_{ij}$ and $K_{ik}$
for $j\neq k$: 
\begin{align*}
\mathbb{V}\left(\hat{f}_{W}(w)\right) & =\mathbb{V}\left(\frac{1}{n}\sum_{i<j}K_{ij}\right)\\
 & =\left(\frac{1}{n}\right)^{2}\sum_{i<j}\sum_{k<l}\mathbb{C}(K_{ij},K_{kl})\\
 & =\left(\frac{1}{n}\right)^{2}[n\cdot\mathcal{\mathbb{C}}(K_{12},K_{12})+2n(N-2)\cdot\mathbb{C}(K_{12},K_{13})]\\
 & =\frac{1}{n}\left[\mathbb{V}(K_{12})+2(N-2)\cdot\mathbb{C}(K_{12},K_{13})\right].
\end{align*}
The third line of this expression uses the fact that, in the summation
in the second line, there are $n=\frac{1}{2}N\left(N-1\right)$ terms
with $(i,j)=(k,l)$ and $N(N-1)(N-2)=2n(N-2)$ terms with one subscript
in common; as noted earlier, when $W_{ij}$ and $W_{kl}$ have no
subscripts in common they are independent (and thus uncorrelated).

To calculate the dependence of this variance on the number of nodes
$N,$ we analyze $\mathbb{V}(K_{12})$ and $\mathbb{C}(K_{12},K_{13}).$
Beginning with the former, 
\begin{align*}
\mathbb{V}(K_{12}) & =\mathbb{E}\left[(K_{12})^{2}\right]-\left(\mathbb{E}[\hat{f}_{W}(w)]\right)^{2}\\
 & =\frac{1}{h^{2}}\int\left[K\left(\frac{w-s}{h}\right)\right]^{2}f_{W}(s)\mathrm{d}s+O(1)\\
 & =\frac{1}{h}\int[K\left(u\right)]^{2}f_{W}(w-hu)\mathrm{d}u+O(1)\\
 & =\frac{f_{W}(w)}{h}\cdot\int[K\left(u\right)]^{2}\mathrm{d}u+O(1)\\
 & \overset{def}{\equiv}\frac{1}{h_{N}}\Omega_{2}(w)+O(1),
\end{align*}
where 
\[
\Omega_{2}(w)\overset{def}{\equiv}f_{W}(w)\cdot\int[K\left(u\right)]^{2}\mathrm{d}u.
\]
Like the expected value, this own variance term is of the same order
of magnitude as in the monadic case, 
\[
\mathbb{V}(K_{12})=O\left(\frac{1}{h}\right).
\]
However, the covariance term $\mathbb{C}(K_{ij},K_{il}),$ which would
be absent for i.i.d. monadic data, is generally nonzero. Since 
\begin{align*}
\mathbb{E}[K_{ij}\cdot K_{ik}] & =\mathbb{E}\left[\int\int\frac{1}{h^{2}}\left[K\left(\frac{w-s_{1}}{h}\right)\right]\cdot\left[K\left(\frac{w-s_{2}}{h}\right)\right]\right.\\
 & \cdot\left.f_{W|AA}(s_{1}|A_{1},A_{2})f_{W|AA}(s_{2}|A_{1},A_{3})\mathrm{d}s_{1}\mathrm{d}s_{2}\right]\\
 & =\mathrm{E}\left[\int\left[K\left(u_{1}\right)\right]f_{W|A}(w-hu_{1}|A_{1})\mathrm{d}u_{1}\right.\\
 & \cdot\left.\int\left[K\left(u_{2}\right)\right]f_{W|A}(w-hu_{2}|A_{1})\mathrm{d}u_{2}\right],\\
 & =\mathbb{E}\left[f_{W|A}(w|A_{1})^{2}\right]+o(1),
\end{align*}
(where the second line uses the change of variables $s_{1}=w-hu_{1}$
and $s_{2}=w-hu_{2}$ and mutual independence of $A_{1},A_{2},$ and
$A_{3}$). It follows that 
\begin{align*}
\mathbb{C}(K_{ij},K_{ik}) & =\mathbb{E}[K_{ij}\cdot K_{ik}]-\left(\mathbb{E}[\hat{f}_{W}(w)]\right)^{2}\\
 & =\left[\mathbb{E}\left[f_{W|A}(w|A_{1})^{2}\right]-f_{W}(w)^{2}\right]+O(h^{2})\\
 & =\mathbb{V}(f_{W|A}(w|A_{1}))+o(1)\\
 & \overset{def}{\equiv}\Omega_{1}(w)+o\left(1\right),
\end{align*}
with 
\[
\Omega_{1}(w)\overset{def}{\equiv}\mathbb{V}(f_{W|A}(w|A_{1})).
\]
Therefore, 
\begin{align}
\mathbb{V}\left(\hat{f}_{W}(w)\right) & =\frac{1}{n}\left[2(N-2)\cdot\mathbb{C}(K_{12},K_{13})+\mathbb{V}(K_{12})\right]\nonumber \\
 & =\frac{4}{N}\Omega_{1}(w)+\left(\frac{1}{nh}\Omega_{2}(w)-\frac{2}{n}\Omega_{1}(w)\right)+o\left(\frac{1}{N}\right)\label{eq: variance}\\
 & =O\left(\frac{4\Omega_{1}(w)}{N}\right)+O\left(\frac{\Omega_{2}(w)}{nh}\right).\nonumber 
\end{align}
and the mean-squared error of $\hat{f}_{W}(w)$ is, using \eqref{eq: bias}
and \eqref{eq: variance}, 
\begin{align}
\mathrm{MSE}\left(\hat{f}_{W}(w)\right)= & \left(\mathbb{E}[\hat{f}_{W}(w)]-f_{W}(w)\right)^{2}+\mathbb{V}\left(\hat{f}_{W}(w)\right)\nonumber \\
= & h^{4}B(w)^{2}+\frac{4}{N}\Omega_{1}(w)+\left(\frac{1}{nh}\Omega_{2}(w)-\frac{2}{n}\Omega_{1}(w)\right)\label{eq: MSE}\\
 & +o(h^{4})+o\left(\frac{1}{N}\right)\nonumber \\
= & O\left(h^{4}\right)+O\left(\frac{4\Omega_{1}(w)}{N}\right)+O\left(\frac{\Omega_{2}(w)}{nh}\right)\nonumber 
\end{align}

Provided that $\Omega_{1}(w)\neq0$ and the bandwidth sequence $h_{N}$
is chosen such that 
\begin{equation}
Nh\rightarrow\infty,\qquad Nh^{4}\rightarrow0\label{eq: MSE_conditions}
\end{equation}
as $N\rightarrow\infty,$ we get that 
\begin{align*}
\mathrm{MSE}\left(\hat{f}_{W}(w)\right) & =o\left(\frac{1}{N}\right)+O\left(\frac{1}{N}\right)+o\left(\frac{1}{N}\right)\\
 & =O\left(\frac{1}{N}\right),
\end{align*}
and hence that 
\[
\sqrt{N}(\hat{f}_{W}(w)-f_{W}(w))=O_{p}(1).
\]
In fact, the rate of convergence of $\hat{f}_{W}(w)$ to $f_{W}(w)$
will be $\sqrt{N}$ as long as $Nh^{4}\leq C\leq Nh$ for some $C>0$
as $N\rightarrow\infty,$ although the mean-squared error will include
an additional bias or variance term of $O(N^{-1})$ if either $Nh$
or $(Nh^{4})^{-1}$ does not diverge to infinity.

To derive the MSE-optimal bandwidth sequence we minimize \eqref{eq: MSE}
with respect to its first and third terms, this yields an optimal
bandwidth sequence of 
\begin{align}
h_{N}^{*}\left(w\right) & =\left[\frac{1}{4}\frac{\Omega_{2}\left(w\right)}{B\left(w\right)^{2}}\frac{1}{n}\right]^{\frac{1}{5}}\label{eq: optimal_bandwidth}\\
 & =O\left(N^{-\frac{2}{5}}\right).\nonumber 
\end{align}
This sequence satisfies condition \eqref{eq: MSE_conditions} above.

Interestingly, the rate of convergence of $\hat{f}_{W}(w)$ to $f_{W}(w)$
under condition \eqref{eq: MSE_conditions} is the same as the rate
of convergence of the sample mean 
\begin{equation}
\bar{W}\overset{def}{\equiv}\frac{1}{n}\sum_{i<j}W_{ij}\label{eq: sample_mean}
\end{equation}
to its expectation $\mu_{W}\overset{def}{\equiv}\mathbb{E}[W_{ij}]$
when $\mathbb{E}[W_{ij}^{2}]<\infty.$ Similar variance calculations
to those for $\hat{f}_{w}(w)$ yield (see also \citet{Holland_Leinhardt_SM76}
and \citet{Menzel_arXiv17}) 
\begin{align*}
\mathbb{V}(\bar{W}) & =O\left(\frac{\mathbb{V}(W_{ij})}{n}\right)+O\left(\frac{4\mathbb{V}(\mathbb{E}[W_{ij}|A_{i}])}{N}\right)\\
 & =O\left(\frac{1}{N}\right),
\end{align*}
provided $\mathbb{E}[W_{ij}|A_{i}]$ is non-degenerate, yielding 
\[
\sqrt{N}(\bar{W}-\mu)=O_{p}(1).
\]
Thus, in contrast to the case of i.i.d monadic data, there is no convergence-rate
``cost'' associated with nonparametric estimation of $f_{W}(w).$
The presence of dyadic dependence, due to its impact on estimation
variance, does slow down the feasible rate of convergence substantially.
With iid data the relevant rate for density estimation would be $n^{2/5}$
when the MSE-optimal bandwidth sequence is used. Recalling that $n=O\left(N^{2}\right)$,
the $\sqrt{N}$ rate we find here corresponds to an $n^{1/4}$ rate.
The slowdown from $n^{2/5}$ to $n^{1/4}$ captures the rate of convergence
costs of dyadic dependence on the variance of our density estimate.

The lack of dependence of the convergence rate of $\hat{f}_{W}(w)$
to $f_{W}(w)$ on the precise bandwidth sequence chosen is analogous
to that for semiparametric estimators defined as averages over nonparametrically-estimated
components \citep[e.g.,][]{Newey_ET94b,Powell_HBE94}. Defining $K_{ji}\overset{def}{\equiv}K_{ij},$
the estimator $\hat{f}_{W}(w)$ can be expressed as 
\[
\hat{f}_{W}(w)=\frac{1}{N}\sum_{i=1}^{N}\hat{f}_{W|A}(w|A_{i}),
\]
where 
\[
\hat{f}_{W|A}(w|A_{i})\overset{def}{\equiv}\frac{1}{N-1}\sum_{j\neq i,j=1}^{N}K_{ij}.
\]
Holding $i$ fixed, the estimator $\hat{f}_{W|A}(W|A_{i})$ can be
shown to converge to $f_{W|A}(w|A_{i})$ at the nonparametric rate
$\sqrt{Nh},$ but the average of this nonparametric estimator over
$A_{i}$ converges at the faster (``parametric'') rate $\sqrt{N}.$
In comparison, while 
\[
\bar{W}=\frac{1}{N}\sum_{i=1}^{N}\hat{\mathbb{E}}\left[\left.W_{ij}\right|A_{i}\right],
\]
for 
\[
\hat{\mathbb{E}}\left[\left.W_{ij}\right|A_{i}\right]\overset{def}{\equiv}\frac{1}{N-1}\sum_{j\neq i,j=1}^{N}W_{ij},
\]
the latter converges at the parametric rate $\sqrt{N},$ and the additional
averaging to obtain $\bar{W}$ does not improve upon that rate.

\section{\label{sec: asymptotic}Asymptotic distribution theory}

To derive conditions under which $\hat{f}_{W}(w)$ is approximately
normally distributed it is helpful to decompose the difference between
$\hat{f}_{W}(w)$ and $f_{W}(w)$ into four terms: 
\begin{align}
\hat{f}_{W}(w)-f_{W}(w) & =\frac{1}{n}\sum_{i<j}(K_{ij}-\mathbb{E}[K_{ij}|A_{i},A_{j}])\label{eq: proj_error1}\\
 & +\frac{1}{n}\sum_{i<j}\mathbb{E}[K_{ij}|A_{i},A_{j}]\label{eq: proj_error2}\\
 & -\left(\mathbb{E}[K_{ij}]+\frac{2}{N}\sum_{i=1}^{N}(\mathbb{E}[K_{ij}|A_{i}]-\mathbb{E}[K_{ij}])\right)\nonumber \\
 & +\frac{2}{N}\sum_{i=1}^{N}(\mathbb{E}[K_{ij}|A_{i}]-\mathbb{E}[K_{ij}])\label{eq: Hajek_proj}\\
 & +\mathbb{E}[K_{ij}]-f_{W}(w)\label{eq: bias_term}\\
 & \equiv T_{1}+T_{2}+T_{3}+T_{4.}\nonumber 
\end{align}
To understand this decomposition observe that the projection of $\hat{f}_{W}(w)=\frac{1}{n}\sum_{i<j}K_{ij}$
onto $\{A_{i}\}_{i=1}^{N}$ equals, by the independence assumptions
imposed on $\{A_{i}\}$ and $\{V_{ij}\},$ the U-statistic $\tbinom{N}{2}^{-1}\sum_{i<j}\mathbb{E}[K_{ij}|A_{i},A_{j}]$.
This U-Statistic is defined in terms of the \emph{latent }i.i.d. random
variables $\{A_{i}\}_{i=1}^{N}$.

The first term in this expression, line \eqref{eq: proj_error1},
is $\hat{f}_{W}(w)$ minus the projection/U-Statistic described above.
Each term in this summation has conditional expectation zero given
the remaining terms (i.e., the terms form a martingale difference
sequence).

The second term in the decomposition, line \eqref{eq: proj_error2},
is the difference between the second-order U-statistic $\frac{1}{n}\sum_{i<j}\mathbb{E}[K_{ij}|A_{i},A_{j}]$
and its\textbf{ }Hájek projection \citep[e.g.,][]{vanderVaart_ASBook00}\footnote{That is the projection of $\frac{1}{n}\sum_{i<j}\mathbb{E}[K_{ij}|A_{i},A_{j}]$
onto the linear subspace consisting of all functions of the form $\sum_{i=1}^{N}g_{i}\left(A_{i}\right)$.}, the third term, line \eqref{eq: Hajek_proj}, is a centered version
of that Hájek projection, and the final term, line \eqref{eq: bias_term},
is the bias of $\hat{f}_{W}(w).$ A similar ``double projection''
argument was used by \citet{Graham_EM17} to analyze the large sample
properties of the Tetrad Logit estimator.

If the bandwidth sequence $h=h_{N}$ satisfies the conditions $Nh\rightarrow\infty$
and $Nh^{4}\rightarrow0,$ the calculations in the previous section
can be used to show that the first, second, and fourth terms of this
decomposition (i.e., $T_{1},$ $T_{2,}$ and $T_{4}$) will all converge
to zero when normalized by $\sqrt{N}$. In this case, $T_{3}$, which
is an average of i.i.d. random variables, will be the leading term
asymptotically such that 
\[
\sqrt{N}(\hat{f}_{W}(w)-f_{W}(w))\overset{D}{\rightarrow}\mathcal{N}(0,4\Omega_{1}(w)),
\]
assuming $\Omega_{1}(w)=\mathbb{V}(f_{W|A}(w|A_{i}))>0$.

If, however, the bandwidth sequence $h$ has $Nh\rightarrow C<\infty$
(a ``knife-edge'' undersmoothing condition similar to one considered
by \citet{Cattaneo_et_al_ET14} in a different context), then both
$T_{1}$ and $T_{3}$ will be asymptotically normal when normalized
by $\sqrt{N}.$ To accommodate both of these cases in a single result,
we will show that a standardized version of the sum $T_{1}+T_{3}$
will have a standard normal limit distribution, although the first,
$T_{1}$, term may be degenerate in the limit.

In Appendix \ref{app: Proofs} we show that both $T_{2}$ and $T_{4}$
will be asymptotically negligible when normalized by the convergence
rate of $T_{1}+T_{3},$ such that the asymptotic distribution of $\hat{f}_{W}(w)$
will only depend on the $T_{1}$ and $T_{3}$ terms.

We start by rewriting the sum of terms $T_{1}$ and $T_{3}$ as 
\begin{align*}
T_{1}+T_{3} & =\frac{1}{n}\sum_{i<j}(K_{ij}-\mathbb{E}[K_{ij}|A_{i},A_{j}])+\frac{2}{N}\sum_{i=1}^{N}(\mathbb{E}[K_{ij}|A_{i}]-\mathbb{E}[K_{ij}])\\
 & \overset{def}{\equiv}\sum_{t=1}^{T(N)}X_{Nt},
\end{align*}
where 
\[
T(N)\equiv N+n
\]
and the triangular array $X_{Nt}$ is defined as 
\begin{align*}
X_{N1} & =\frac{2}{N}\left(\mathbb{E}[K_{12}|A_{1}]-\mathbb{E}[K_{12}]\right),\\
X_{N2} & =\frac{2}{N}\left(\mathbb{E}[K_{23}|A_{2}]-\mathbb{E}[K_{23}]\right),\\
 & \vdots\\
X_{NN} & =\frac{2}{N}(E[K_{N,1}|A_{N}]-\mathbb{E}[K_{N,1}]),\\
X_{N,N+1} & =\frac{1}{n}(K_{12}-\mathbb{E}[K_{12}|A_{1},A_{2}]),\\
X_{N,N+2} & =\frac{1}{n}(K_{13}-\mathbb{E}[K_{13}|A_{1},A_{3}])\\
 & \vdots\\
X_{N,N+N-1} & =\frac{1}{n}(K_{1N}-\mathbb{E}[K_{1N}|A_{1},A_{N}]),\\
 & \vdots\\
X_{N,N+n} & =\frac{1}{n}(K_{N-1,N}-\mathbb{E}[K_{N-1,N}|A_{N-1},A_{N}]).
\end{align*}
That is, $\{X_{Nt}\}$ is the collection of terms of the form 
\[
\frac{2}{N}(\mathbb{E}[K_{ij}|A_{i}]-\mathbb{E}[K_{ij}])
\]
for $i=1,...,N$ (with $j\neq i$) and 
\[
\frac{1}{n}(K_{ij}-\mathbb{E}[K_{ij}|A_{i},A_{j}])
\]
for $i=1,...,N-1$ and $j=i+1,...,N.$ Using the independence assumptions
on $\{A_{i}\}_{i=1}^{N}$ and $\{V_{ij}\}_{i<j}$, as well as iterated
expectations, it is tedious but straightforward to verify that 
\[
\mathbb{E}[X_{Nt}|\{X_{Ns},s\neq t\}]=0,
\]
that is, $X_{NT}$ is a martingale difference sequence (MDS).

Defining the variance of this MDS as 
\begin{align*}
\sigma_{N}^{2} & \overset{def}{\equiv}\mathbb{E}\left(\sum_{t=1}^{T(N)}X_{Nt}\right)^{2}\\
 & =\sum_{t=1}^{T(N)}\mathbb{V}(X_{Nt}),
\end{align*}
we can demonstrate asymptotic normality of its standardized sum --
$\frac{1}{\sigma_{N}}\sum_{t=1}^{T(N)}X_{Nt}$ -- by a central limit
theorem for martingale difference triangular arrays (see, for example,
\citet{Hall_Heyde_Bk1980}, Theorem 3.2 and Corollary 3.1 and \citet{White_Bk01},
Theorem 5.24 and Corollary 5.26). Specifically, if the Lyapunov condition
\begin{equation}
\sum_{t=1}^{T(N)}\mathbb{E}\left(\frac{X_{Nt}}{\sigma_{N}}\right)^{r}\rightarrow0\label{eq: Liapunov}
\end{equation}
holds for some $r>2,$ and also the stability condition
\begin{equation}
\sum_{t=1}^{T(N)}\left(\frac{X_{Nt}}{\sigma_{N}}\right)^{2}\overset{p}{\rightarrow}1,\label{eq: mixing}
\end{equation}
holds then 
\begin{align}
\sum_{t=1}^{T(N)}\frac{X_{Nt}}{\sigma_{N}} & =\frac{1}{\sigma_{N}}(T_{1}+T_{3})\nonumber \\
 & \overset{D}{\rightarrow}\mathcal{N}(0,1).\label{eq: MDCLT}
\end{align}
From the calculations used in the MSE analysis of Section \ref{sec: MSE}
we have that 
\begin{align*}
\sigma_{N}^{2} & =\mathbb{V}(T_{1})+\mathbb{V}(T_{3})\\
 & =\frac{\mathbb{E}[K_{ij}^{2}]}{n}+\frac{4\mathbb{V}(\mathbb{E}[K_{ij}|A_{i}])}{N}+O\left(\frac{1}{n}\right)\\
 & =\frac{\Omega_{2}(w)}{nh}+\frac{4\Omega_{1}(w)}{N}+O\left(\frac{1}{n}\right)+O\left(\frac{h^{2}}{N}\right),
\end{align*}
so, taking $r=3$, 
\[
\frac{1}{\sigma_{N}^{2}}=O(N)
\]
assuming $\Omega_{1}(w)>0$ and $Nh\geq C>0.$ In the degenerate case,
where $\mathbb{V}(\mathbb{E}[K_{ij}|A_{i}])=\Omega_{1}(w)=0$, we
will still have $(\sigma_{N})^{-2}=O(nh)=O(N)$ as long as the ``knife-edge''
$h\propto N^{-1}$ undersmoothing bandwidth sequence is chosen.

To verify the Lyapunov condition (\ref{eq: Liapunov}), note that
\begin{align}
\mathbb{E}\left(\frac{1}{n}(K_{ij}-\mathbb{E}[K_{ij}|A_{i},A_{j}])\right)^{3} & \leq8\mathbb{E}\left(\frac{K_{ij}}{n}\right)^{3}\nonumber \\
 & =\frac{8}{n^{3}}\frac{1}{h^{3}}\int\left[K\left(\frac{w-s}{h}\right)\right]^{3}f_{W}(s)ds\nonumber \\
 & =\frac{8}{n^{3}h^{2}}\int[K\left(u\right)]^{3}f_{W}(w-hu)du\nonumber \\
 & =O\left(\frac{1}{n^{3}h^{2}}\right)\label{eq: third moment 1}
\end{align}
and 
\begin{align}
\mathbb{E}\left(\frac{2}{N}(\mathbb{E}[K_{ij}|A_{i}]-\mathbb{E}[K_{ij}])\right)^{3} & \leq\frac{8^{2}}{N^{3}}\mathbb{E}\left(\mathbb{E}[K_{ij}|A_{i}]\right)^{3}\nonumber \\
 & =\frac{8^{2}}{N^{3}}\mathbb{E}\left(\int K\left(u\right)f_{W|A}(w-hu|A_{i})du\right)^{3}\nonumber \\
 & =O\left(\frac{1}{N^{3}}\right).\label{eq: third moment 2}
\end{align}
Putting things together we get that 
\begin{align*}
\sum_{t=1}^{T(N)}\mathbb{E}\left(X_{Nt}\right)^{3}= & n\mathbb{E}\left(\frac{1}{n}(K_{ij}-\mathbb{E}[K_{ij}|A_{i},A_{j}])\right)^{3}\\
 & +N\mathbb{E}\left(\frac{2}{N}(\mathbb{E}[K_{ij}|A_{i}]-\mathbb{E}[K_{ij}])\right)^{3}\\
= & O\left(\frac{1}{\left(nh\right)^{2}}\right)+O\left(\frac{1}{N^{2}}\right)\\
= & O\left(\frac{1}{N^{2}}\right)
\end{align*}
when $Nh\geq C>0$ for all $N.$ Therefore the Lyapunov condition
(\ref{eq: Liapunov}) is satisfied for $r=3,$ since
\begin{align*}
\sum_{t=1}^{T(N)}\mathbb{E}\left(\frac{X_{Nt}}{\sigma_{N}}\right)^{3} & =O(N^{3/2})\cdot O\left(\frac{1}{N^{2}}\right)\\
 & =O\left(\frac{1}{\sqrt{N}}\right)\\
 & =o(1).
\end{align*}
To verify the stability condition (\ref{eq: mixing}), we first rewrite
that condition as
\begin{align*}
0 & =\underset{N\rightarrow\infty}{\lim}\left(\frac{1}{\sigma_{N}^{2}}\sum_{t=1}^{T\left(N\right)}\left(X_{Nt}^{2}-\mathbb{E}\left[X_{Nt}^{2}\right]\right)\right)\\
 & =\underset{N\rightarrow\infty}{\lim}\left(\frac{1}{N\sigma_{N}^{2}}\sum_{t=1}^{T\left(N\right)}\left(R_{1}+R_{2}\right)\right)
\end{align*}
where
\begin{align*}
R_{1}\equiv & N\sum_{i=1}^{N}\left[\left(\frac{2}{N}(\mathbb{E}[K_{ij}|A_{i}]-\mathbb{E}[K_{ij}])\right)^{2}-\mathbb{E}\left(\frac{2}{N}(\mathbb{E}[K_{ij}|A_{i}]-\mathbb{E}[K_{ij}])\right)^{2}\right]\\
 & +N\sum_{i<j}\left[\left(\frac{1}{n}(K_{ij}-\mathbb{E}[K_{ij}|A_{i},A_{j}])\right)^{2}-\mathbb{E}\left[\left.\left(\frac{1}{n}(K_{ij}-\mathbb{E}[K_{ij}|A_{i},A_{j}])\right)^{2}\right\vert A_{i},A_{j}\right]\right]
\end{align*}
and
\[
R_{2}\equiv N\sum_{i<j}\left[\mathbb{E}\left[\left.\left(\frac{1}{n}(K_{ij}-\mathbb{E}[K_{ij}|A_{i},A_{j}])\right)^{2}\right\vert A_{i},A_{j}\right]-\mathbb{E}\left[\left(\frac{1}{n}(K_{ij}-\mathbb{E}[K_{ij}|A_{i},A_{j}])\right)^{2}\right]\right].
\]
Since $1/N\sigma_{N}^{2}=O(1),$ the stability condition (\ref{eq: mixing 2})
will hold if $R_{1}$ and $R_{2}$ both converge to zero in probability.

By the independence restrictions on $\{U_{ij}\}$ and $\{A_{i}\},$
the (mean zero) summands in $R_{1}$ are mutually uncorrelated, so
\begin{align*}
\mathbb{E}\left[R_{1}^{2}\right] & \equiv N^{2}\sum_{i=1}^{N}\mathbb{E}\left[\left(\left(\frac{2}{N}(\mathbb{E}[K_{ij}|A_{i}]-\mathbb{E}[K_{ij}])\right)^{2}-\mathbb{E}\left(\frac{2}{N}(\mathbb{E}[K_{ij}|A_{i}]-\mathbb{E}[K_{ij}])\right)^{2}\right)^{2}\right]\\
 & +N^{2}\sum_{i<j}\mathbb{E}\left[\left(\left(\frac{1}{n}(K_{ij}-\mathbb{E}[K_{ij}|A_{i},A_{j}])\right)^{2}-\mathbb{E}\left[\left.\left(\frac{1}{n}(K_{ij}-\mathbb{E}[K_{ij}|A_{i},A_{j}])\right)^{2}\right\vert A_{i},A_{j}\right]\right)^{2}\right]\\
 & =O\left(\frac{\mathbb{E}\left(\mathbb{E}[K_{ij}|A_{i}]\right)^{4}}{N}\right)+O\left(\frac{N^{2}\mathbb{E}\left(K_{ij}\right)^{4}}{n^{3}}\right).
\end{align*}
But, using analogous arguments to (\ref{eq: third moment 1}) and
((\ref{eq: third moment 2}), 
\[
\mathbb{E}\left[\mathbb{E}[K_{ij}|A_{i}]^{4}\right]=O\left(1\right)
\]
and
\[
\mathbb{E}\left[K_{ij}^{4}\right]=O\left(\frac{1}{h^{3}}\right),
\]
so
\begin{align*}
\mathbb{E}\left[R_{1}^{2}\right] & =O\left(\frac{1}{N}\right)+O\left(\frac{N^{2}}{(nh)^{3}}\right)\\
 & =O\left(\frac{1}{N}\right)\\
 & =o(1),
\end{align*}
under the bandwidth condition that $1/nh=O(1/N).$ So $R_{1}$ converges
in probability to zero. Moreover, $R_{2}$ is proportional to a (mean
zero) second-order U-statistic,
\begin{align*}
R_{2} & =\frac{1}{n}\sum_{i<j}\frac{N}{n}\left[\mathbb{E}\left[\left.(K_{ij}-\mathbb{E}[K_{ij}|A_{i},A_{j}])^{2}\right\vert A_{i},A_{j}\right]-\mathbb{E}\left[(K_{ij}-\mathbb{E}[K_{ij}|A_{i},A_{j}])^{2}\right]\right]\\
 & \equiv\frac{1}{n}\sum_{i<j}p_{N}(A_{i},A_{j}),
\end{align*}
with kernel having second moment
\begin{align*}
\mathbb{E}\left[p_{N}(A_{i},A_{j})^{2}\right] & =O\left(\frac{N^{2}}{n^{2}}\mathbb{E}\left(\mathbb{E}[K_{ij}^{2}|A_{i},A_{j}]\right)^{2}\right)\\
 & =O\left(\frac{N^{2}}{n^{2}}\mathbb{\cdot}\frac{1}{h^{2}}\right)\\
 & =O(1)\\
 & =o(N),
\end{align*}
again imposing the bandwidth restriction $1/nh=O(1/N)$. Thus by Lemma
3.1 of \citet{Powell_Stock_Stoker_EM89}, $R_{2}$ converges in probability
to its (zero) expected value.

Since conditions (\ref{eq: Liapunov}) and (\ref{eq: mixing}) both
hold, a central limit theorem for martingale difference triangular
arrays implies
\[
\frac{1}{\sigma_{N}}(T_{1}+T_{3})\overset{D}{\rightarrow}\mathcal{N}(0,1).
\]
A final step is to used this result to obtain the asymptotic distribution
of $\hat{f}_{W}(w).$ Because 
\[
\frac{1}{\sigma_{N}}=O\left(\sqrt{N}\right),
\]
we have that $T_{2}$ and $T_{4}$ are asymptotically negligible after
standardization with $\sigma_{N}^{-1}$ (see Appendix \ref{app: Proofs}),
\[
\frac{T_{2}}{\sigma_{N}}=O_{p}\left(\sqrt{\frac{N}{n}}\right)=o_{p}(1)
\]
and 
\[
\frac{T_{4}}{\sigma_{N}}=O\left(\sqrt{N}h^{2}\right)=o(1),
\]
so that 
\begin{align*}
\frac{1}{\sigma_{N}}\left(\hat{f}_{W}(w)-f_{W}(w)\right) & =\frac{1}{\sigma_{N}}\left(T_{1}+T_{2}+T_{3}+T_{4}\right)\\
 & \overset{D}{\rightarrow}\mathcal{N}(0,1).
\end{align*}
When $Nh^{4}\rightarrow0$ and $Nh\rightarrow\infty,$ 
\[
N\sigma_{N}^{2}\rightarrow4\Omega_{1}(w)
\]
and 
\[
\sqrt{N}\left(\hat{f}_{W}(w)-f_{W}(w)\right)\overset{D}{\rightarrow}\mathcal{N}(0,4\Omega_{1}(w))
\]
as long as $\mathbb{V}(\mathbb{E}[K_{ij}|A_{i}])>0.$

Under ``knife-edge'' bandwidth sequences, such that $Nh\rightarrow C>0,$
we have instead that 
\[
N\sigma_{N}^{2}\rightarrow4\Omega_{1}(w)+C^{-1}\Omega_{2}(w)
\]
and 
\[
\sqrt{N}(\hat{f}_{W}(w)-f_{W}(w))\overset{D}{\rightarrow}\mathcal{N}(0,4\Omega_{1}(w)+C^{-1}\Omega_{2}(w)).
\]

\subsubsection*{Degeneracy}

Degeneracy arises when $\mathbb{V}(\mathbb{E}[K_{ij}|A_{i}])=\Omega_{1}\left(w\right)=0.$
In terms of the underlying network generating process (NGP), degeneracy
arises when the conditional density of $W_{ij}$ at $w$ given $A_{i}=a$
is constant in $a$ (i.e., when $\mathbb{V}\left(f_{W|A}\left(\left.w\right|A_{i}\right)\right)=0$).

As a simple example of such an NGP, let $A_{i}$ equal $-1$ with
probability $\pi$ and $1$ otherwise; next set 
\[
W_{ij}=A_{i}A_{j}+V_{ij}
\]
with $V_{ij}$ standard normal. In this case the conditional density
$f_{W|A}\left(\left.w\right|A_{i}\right)$ is the mixture 
\[
f_{W|A}\left(\left.w\right|A_{i}\right)=\pi\phi\left(w+A_{i}\right)+\left(1-\pi\right)\phi\left(w-A_{i}\right)
\]
with $\phi\left(\cdot\right)$ the standard normal density function.
Unconditionally the density is 
\[
f_{W}\left(w\right)=\left[\pi^{2}+\left(1-\pi\right)^{2}\right]\phi\left(w-1\right)+2\pi\left(1-\pi\right)\phi\left(w+1\right).
\]
Observe that, if $\pi=1/2$, then $f_{W|A}\left(\left.w\right|A_{i}=1\right)=f_{W|A}\left(\left.w\right|A_{i}=-1\right)=f_{W}\left(w\right)$
and hence that $\mathbb{V}\left(f_{W|A}\left(\left.w\right|A_{i}\right)\right)=0$.\footnote{Degeneracy also arises when $w=1$.}
Degeneracy arises in this case, even though there is non-trivial dependence
across dyads sharing an agent in common. If $\pi\neq1/2$, then $\mathbb{V}\left(f_{W|A}\left(\left.w\right|A_{i}\right)\right)>0$,
but one still might worry about ``near degeneracy'' when $\pi$
is close to $1/2$.

\citet{Menzel_arXiv17} shows that under degeneracy, the limit distribution
of the sample mean, $\bar{W}$, equation \eqref{eq: sample_mean}
on \vpageref{eq: sample_mean} above, may be non-Gaussian. This occurs
because (i) the $T_{1}$ and $T_{2}$ terms in a double projection
decomposition of $\bar{W}$ analogous to the one used here for $\hat{f}_{W}\left(w\right)$
will be of equal order \emph{and} $T_{2}$, the Hájek Projection error,
may be non-Gaussian (as is familiar from the theory of U-Statistics,
e.g., Chapter 12 of \citet{vanderVaart_ASBook00}).

The situation is both more complicated and simpler here. In the case
of the estimated density $\hat{f}_{W}\left(w\right)$, if the bandwidth
sequence $h=h_{N}$ satisfies the conditions $Nh\rightarrow\infty$
and $Nh^{4}\rightarrow0,$ then $T_{2}$ will be of smaller order
than $T_{1}$ and hence not contribute to the limit distribution irrespective
of whether the NGP is degenerate or not. In particular, under degeneracy
the Liaponuv condition \eqref{eq: Liapunov} continues to hold for
$r=3$ since 
\[
\sum_{t=1}^{T(N)}E\left(\frac{X_{Nt}}{\sigma_{N}}\right)^{3}=O\left(\frac{1}{\sqrt{nh}}\right)
\]
and it follows straightforwardly that $\frac{1}{\sigma_{N}}\left(\hat{f}_{W}\left(w\right)-f_{W}\left(w\right)\right)$
continues to be normal in the limit.

The ``knife-edge'' undersmoothing bandwidth sequence is primarily
of interest because it results in a sequence where both $T_{1}$ and
$T_{3}$ contribute to the limit distribution. In practice this does
not mean that the researcher should set $h=h_{N}\propto N^{-1}$.
Based on the theoretical analysis sketched above, we recommend choosing
a sequence that tends to zero slightly faster than mean squared error
optimal sequence where $h=h_{N}\propto n^{-1/5}$.\footnote{In practice ``plug-in'' bandwidths that would be appropriate in
the absence of any dyadic dependence across the $\left\{ W_{ij}\right\} _{i<j}$
might work well; although this remains an unexplored conjecture.}

Under such a sequence we will have 
\[
\sqrt{N}(\hat{f}_{W}(w)-f_{W}(w))\overset{D}{\rightarrow}\mathcal{N}(0,4\Omega_{1}(w))
\]
under non-degeneracy and 
\[
\sqrt{nh}(\hat{f}_{W}(w)-f_{W}(w))\overset{D}{\rightarrow}\mathcal{N}(0,\Omega_{2}(w))
\]
under degeneracy. Although the rate of convergence of $\hat{f}_{W}(w)$
to $f_{W}(w)$ is faster in the case of degeneracy this will not affect
inference in practice as long as an appropriate estimate of $\sigma_{N}$
is used; that is working directly with $(\hat{f}_{W}(w)-f_{W}(w))/\sigma_{N}$
ensures rate-adaptivity. Note also that, in the absence of degeneracy,
the MSE optimal bandwidth sequence could be used. By slightly undersmoothing
relative to this sequence, we ensure that the limit distribution remains
unbiased in case of degeneracy.

\section{\label{sec: asymptotic_variance_estimation}Asymptotic variance estimation}

To construct Wald-based confidence intervals for $\hat{f}_{W}(w),$
a consistent estimator of its asymptotic variance is needed. When
$Nh\rightarrow C<\infty,$ the asymptotic variance depends on both
\[
\Omega_{2}(w)\overset{def}{\equiv}f_{W}(w)\cdot\int[K\left(u\right)]^{2}\mathrm{d}u
\]
and 
\[
\Omega_{1}(w)\overset{def}{\equiv}\mathbb{V}\left(f_{W|A}(w|A_{i})\right).
\]
In this section we present consistent estimators for both of these
terms.

A simple estimator of $\Omega_{2}(w)$ is 
\begin{align}
\tilde{\Omega}_{2}(w) & =\frac{h}{n}\sum_{i<j}K_{ij}^{2},\label{eq: OMEGA2_tilde}
\end{align}
the consistency of which we demonstrate in Appendix \ref{app: Proofs}:

\begin{equation}
\tilde{\Omega}_{2}(w)\overset{p}{\rightarrow}\Omega_{2}(w).\label{eq: consistency_of_OMEGA2_tilde}
\end{equation}
The estimator $\tilde{\Omega}_{2}(w)$ uses the second moment of $K_{ij}$
instead of its sample variance to estimate $\Omega_{2}(w);$ in practice
we recommend, similar to \citet{Newey_ET94b} in the context of monadic
kernel-based estimation, the less conservative alternative: 
\begin{align*}
\hat{\Omega}_{2}(w) & \equiv h\left(\left(\frac{1}{n}\sum_{i<j}K_{ij}^{2}\right)-\left(\hat{f}_{W}(w)\right)^{2}\right)\\
 & =h\left(\frac{1}{n}\sum_{i<j}\left(K_{ij}-\hat{f}_{W}(w)\right)^{2}\right)\\
 & =\tilde{\Omega}_{1}(w)+o_{p}(1)\\
 & =\Omega_{1}(w)+o_{p}(1).
\end{align*}

We next turn to estimation of 
\[
\Omega_{1}(w)=\mathbb{V}\left(f_{W|A}(w|A_{1})\right)=\lim_{N\rightarrow\infty}\mathbb{C}(K_{ij},K_{ij})
\]
where $i\neq k.$ A natural sample analog estimator, following a suggestion
by \citet{Graham_HBE18} in the context of parametric dyadic regression,
involves an average over the three indices $i,$ $j,$ and $k$: 
\begin{align*}
\hat{\Omega}_{1}(w) & \equiv\frac{1}{N(N-1)(N-2)}\sum_{i\neq j\neq k}(K_{ij}-\hat{f}_{W}(w))(K_{ik}-\hat{f}_{W}(w))\\
 & \equiv\left(\begin{array}{c}
N\\
3
\end{array}\right)^{-1}\sum_{i<j<k}S_{ijk}-\hat{f}_{W}(w)^{2},
\end{align*}
for $S_{ijk}=\frac{1}{3}\left(K_{ij}K_{ik}+K_{ij}K_{jk}+K_{ik}K_{jk}\right).$\footnote{See also the variance estimator for density presented in \citet{Holland_Leinhardt_SM76}.}
In Appendix \ref{app: Proofs} we show that 
\begin{equation}
\hat{\Omega}_{1}(w)\overset{p}{\rightarrow}\Omega_{1}(w).\label{eq: consistency_of_OMEGA1_hat}
\end{equation}

Inserting these estimators, $\hat{\Omega}_{1}(w)$ and $\hat{\Omega}_{2}(w)$,
into the formula for the variance of $\hat{f}_{W}(w)$ yields a variance
estimate of 
\begin{equation}
\hat{\sigma}_{N}^{2}=\frac{1}{nh}\hat{\Omega}_{2}(w)+\frac{2(N-2)}{n}\hat{\Omega}_{1}(w).\label{eq: var_f(w)_hat}
\end{equation}
We end this section by observing that the following equality holds
\begin{align*}
\hat{\sigma}_{N}^{2}= & \frac{1}{n^{2}}\sum_{i<j}\left(K_{ij}-\hat{f}_{W}(w)\right)^{2}\\
 & +\frac{2(N-2)}{n}\left(\frac{1}{N(N-1)(N-2)}\sum_{i\neq j\neq k}(K_{ij}-\hat{f}_{W}(w))(K_{ik}-\hat{f}_{W}(w))\right)\\
= & \frac{1}{n^{2}}\left(\sum_{i<j}\sum_{k<l}d_{ijkl}(K_{ij}-\hat{f}_{W}(w))(K_{kl}-\hat{f}_{W}(w))\right),
\end{align*}
where 
\[
d_{ijkl}=1\{i=j,k=l,i=l,\text{or }j=k\}.
\]
As \citet{Graham_HBE18} notes, this coincides with the estimator
for 
\[
\mathbb{V}(\bar{W})=\mathbb{V}\left(\frac{1}{n}\sum_{i<j}W_{ij}\right)
\]
proposed by \citet{Fafchamp_Gubert_JDE07}, replacing ``$W_{ij}-\bar{W}$''
with ``$K_{ij}-\bar{K}$'', with $\bar{K}\overset{def}{\equiv}\hat{f}_{W}(w)$
(see also \citet{Holland_Leinhardt_SM76}, \citet{Cameron_Miller_WP14}
and \citet{Aronow_et_al_PA17}). Our variance estimator can also be
viewed as a dyadic generalization of the variance estimate proposed
by \citet{Newey_ET94b} for ``monadic'' kernel estimates.

\section{\label{sec: simulation_study}Simulation study}

Our simulations design is based upon the example used to discuss degeneracy
in Section \ref{sec: asymptotic}. As there we let $A_{i}$ equal
$-1$ with probability $\pi$ and $1$ otherwise. We generate $W_{ij}$
\[
W_{ij}=A_{i}A_{j}+V_{ij}
\]
with $V_{ij}\sim\mathcal{N}(0,1)$. We set $\pi=1/3$ and estimate
the density $f_{W}\left(w\right)$ at $w=1.645$.

We present results for three sample sizes: $N=100,400$ and $1,600$.
These sample sizes are such that, for a ``sufficiently non-degenerate''
NGP, the standard error of $\hat{f}_{W}\left(w\right)$ would be expected
to decline by a factor of $1/2$ for each increase in sample size
(if the bandwidth is large enough to ensure that the $\frac{\Omega_{2}\left(w\right)}{nh}$
variance term is negligible relative to the $\frac{2\Omega_{1}\left(w\right)\left(N-2\right)}{n}\approx\frac{4\Omega_{1}\left(w\right)}{N}$
one). We set the bandwidth equal to the MSE-optimal one presented
in equation \eqref{eq: optimal_bandwidth} above. This is an `oracle'
bandwidth choice. Developing feasible data-based methods of bandwidth
selection would be an interesting topic for future research.

Table \ref{tab: monte_carlo_designs} presents the main elements of
each simulation design. Panel B of the table lists ``pencil and paper''
bias and asymptotic standard error calculations based upon the expressions
presented in Section \ref{sec: MSE} above. Panel B also presents
analytic estimates of the standard deviations of the $T_{1}$ and
$T_{3}$ terms in the decomposition of $\hat{f}_{W}\left(w\right)$
used to derive its limit distribution. In the given designs both terms
of are similar magnitude despite the fact that the contribution of
the $T_{1}$ term is asymptotically negligible in theory.

\begin{table}
\caption{Monte Carlo Designs}
\label{tab: monte_carlo_designs}
\begin{centering}
\begin{tabular}{|c|c|c|c|}
\hline 
$N$  & 100  & 400  & 1,600\tabularnewline
\hline 
\hline 
\multicolumn{4}{|c|}{\textbf{Panel A: Design \& Bandwidth}}\tabularnewline
\hline 
$\pi$  & $\frac{1}{3}$  & $\frac{1}{3}$  & $\frac{1}{3}$\tabularnewline
\hline 
$w$  & 1.645  & 1.645  & 1.645\tabularnewline
\hline 
$h_{N}^{*}\left(w\right)$  & 0.2496  & 0.1431  & 0.0822\tabularnewline
\hline 
\multicolumn{4}{|c|}{\textbf{Panel B: Theoretical Sampling Properties}}\tabularnewline
\hline 
$h^{2}B(w)$  & -0.0033  & -0.0011  & -0.0004\tabularnewline
\hline 
$\mathrm{ase}\left(\hat{f}_{W}\left(w\right)\right)=\sqrt{\frac{2\Omega_{1}\left(w\right)\left(N-2\right)}{n}+\frac{\Omega_{2}\left(w\right)}{nh}}$  & 0.0117  & 0.0053  & 0.0025\tabularnewline
\hline 
$\mathrm{ase}\left(T_{3}\right)=\sqrt{\frac{2\Omega_{1}\left(w\right)\left(N-2\right)}{n}}$  & 0.0098  & 0.0049  & 0.0025\tabularnewline
\hline 
$\mathrm{ase}\left(T_{1}\right)=\sqrt{\frac{\Omega_{2}\left(w\right)}{nh}}$  & 0.0065  & 0.0021  & 0.0007\tabularnewline
\hline 
\end{tabular}
\par\end{centering}
\textsc{\uline{Notes}}: Rows $1$ through $3$ list the basic Monte
Carlo design and bandwidth parameter choices. The bandwidths coincide
with the MSE optimal one given in equation \eqref{eq: optimal_bandwidth}.
Panel B gives pencil and paper calculations for the bias of $\hat{f}_{W}\left(w\right)$,
as well as its asymptotic standard error (ase), based upon, respectively,
equations \eqref{eq: bias} and \eqref{eq: variance} in Section \ref{sec: MSE}.
The asymptotic standard errors of $T_{1}$ and $T_{3}$, as defined
in Section \ref{sec: asymptotic}, are also separately given. 
\end{table}

Table \ref{tab: monte_carlo_results} summarizes the results of 1,000
Monte Carlo simulations. The median bias and standard deviation of
our density estimates across the Monte Carlo replications closely
track our theoretical predictions (compare rows 1 and 2 of Table \ref{tab: monte_carlo_results}
with Rows 1 and 2 of Panel B of Table \ref{tab: monte_carlo_designs}.
Row 3 of the table reports the median ``Fafchamps and Gubert'' asymptotic
standard error estimate. This standard error estimate is generally
larger than its asymptotic counterpart. Consequently the coverage
of confidence intervals based upon it is conservative (Row 5). The
degree of conservatism is declining in sample size, suggesting that
-- as expected -- the ``Fafchamps and Gubert'' asymptotic standard
error estimate is closer to its theoretical counterpart as $N$ grows.
Row 4 of the table reports the coverage of confidence intervals based
upon standard errors which ignore the presence of dyadic dependence;
these intervals -- as expected -- fail to cover the true density
frequently enough.

\begin{table}
\caption{Monte Carlo Results}
\label{tab: monte_carlo_results}
\begin{centering}
\begin{tabular}{|c|c|c|c|}
\hline 
$N$  & $100$  & $400$  & $1,600$\tabularnewline
\hline 
\hline 
median bias  & -0.0028  & -0.0010  & -0.0006\tabularnewline
\hline 
standard deviation  & 0.0112  & 0.0051  & 0.0025\tabularnewline
\hline 
median $\mathrm{\hat{ase}}\left(\hat{f}_{W}\left(w\right)\right)$  & 0.0173  & 0.0068  & 0.0028\tabularnewline
\hline 
coverage (iid)  & 0.678  & 0.551  & 0.390\tabularnewline
\hline 
coverage (FG)  & 0.995  & 0.987  & 0.967\tabularnewline
\hline 
\end{tabular}
\par\end{centering}
\textsc{\uline{Notes}}: A robust measure of the standard deviation
of $\hat{f}_{W}\left(w\right)$ is reported in row 2. It equals the
difference between the 0.95 and 0.05 quantiles of the Monte Carlo
distribution of $\hat{f}_{W}\left(w\right)$ divided by $2\times1.645$.
Row 4 reports the coverage of a nominal 95 percent Wald-based confidence
interval that ignores the presence of dyadic dependence. Row 5 reports
the coverage properties of a nominal 95 percent Wald-based confidence
interval that uses the \citet{Fafchamp_Gubert_JDE07} variance estimate
discussed in Section \ref{sec: asymptotic}. 
\end{table}

The simulations suggest, for the designs considered, that the asymptotic
theory presented in Sections \ref{sec: MSE} and \ref{sec: asymptotic}
provides an accurate approximation of finite sample behavior. Our
variance estimate is a bit conservative for the designs considered;
whether this is peculiar to the specific design considered or a generic
feature of the estimate is unknown.\footnote{We observe that our variance estimate implicitly includes an estimate
of the variance of $T_{2}$, which is negligible in the limit.} As with bandwidth selection, further exploration of methods of variance
estimation in the presence of dyadic dependence is warranted.

\section{\label{sec: extensions}Extensions}

There are a number of avenues for extension or modification of the
simple results for scalar density estimation presented above. One
variant of these results would apply when the dyadic variable $W_{ij}$
lacks the idiosyncratic component $V_{ij},$ i.e., when 
\[
W_{ij}=W(A_{i},A_{j}),
\]
for $\{A_{i}\}$ an i.i.d. sequence. This case arises when $W_{ij}$
is a measure of ``distance'' between the attributes of nodes $i$
and $j,$ for example, 
\[
W_{ij}=\sqrt{\left(A_{i}-A_{j}\right)^{2}},
\]
for $A_{i}$ a scalar measure of ``location'' for agent $i.$ The
asymptotic distribution of $\hat{f}_{W}(w)$ derived above should
be applicable to this case as long as the conditional density function
$f_{W|A}(w|a)$ of $W_{ij}$ given $A_{i}$ is well-defined, which
would be implied if $A_{i}$ has a continuously-distributed component
given its remaining component (if any) and the function $W(\cdot)$
is continuously differentiable in that component. In the decomposition
of $\hat{f}_{W}(w)-f_{W}(w)$ for this case, the term corresponding
to $T_{1}$ would be identically zero (as would $\Omega_{2}(w)$),
but the $T_{2}$ term could still be shown to be asymptotically negligible
using Lemma 3.1 of \citet{Powell_Stock_Stoker_EM89} as long as $Nh\rightarrow\infty.$

Another straightforward extension of this analysis would be to directed
dyadic data, where $W_{ij}$ is observed for all pairs of indices
with $i\neq j$ and $W_{ij}\neq W_{ji}$ with positive probability.
The natural generalization of the data generation process would be
\[
W_{ij}=W(A_{i},B_{j},V_{ij}),
\]
with $\{A_{i}\},$ $\{B_{j}\},$ and $\{V_{ij}\}$ mutually independent
and i.i.d. with $V_{ij}\neq V_{ji}$ in general. Here the conditional
densities 
\[
f_{W|A}(w|a)=\mathbb{E}[f_{W|AB}(w|A_{i}=a,B_{j})]
\]
and 
\[
f_{W|B}(w|b)=\mathbb{E}[f_{W|AB}(w|A_{i},B_{j}=b)]
\]
will differ, and the asymptotic variance of $\hat{f}_{W}(w)$ will
depend upon 
\[
\Omega_{1}(w)=\mathbb{V}\left(\frac{1}{2}\left(f_{W|A}(w|A_{i})+f_{W|B}(w|B_{i})\right)\right)
\]
in a way analogous to how $\Omega_{1}(w)$, defined earlier, does
in the undirected case analyzed in this paper.

Yet another generalization of the results would allow $W_{ij}$ to
be a $p$-dimensional jointly-continuous $W_{ij}$ random vector.
The estimator 
\[
\hat{f}_{W}(w)=\frac{1}{n}\sum_{i=1}^{N-1}\sum_{j=1+1}^{N}\frac{1}{h^{p}}K\left(\frac{w-W_{ij}}{h}\right)
\]
of the $p$-dimensional density function $f_{W}(w)$ will continue
to have the same form as derived in the scalar case, provided $Nh^{p}\rightarrow\infty$
(or $Nh^{p}\rightarrow C>0$) as long as the relevant bias term $T_{4}$
is negligible. If the density is sufficiently smooth and $K(\cdot)$
is a "higher-order kernel" with 
\begin{align*}
\int K(u)\mathrm{d}u & =1,\\
\int u_{1}^{j_{1}}u_{2}^{j_{2}}...u_{p}^{j_{p}}K(u)\mathrm{d}u & =0\qquad\text{for }j_{i}\in\{0,...,q\}\text{ with }\sum_{i=1}^{p}j_{i}<q,
\end{align*}
then the bias term $T_{4}$ will satisfy 
\begin{align*}
T_{4} & \equiv\mathbb{E}\left[\hat{f}_{W}(w)\right]-f_{W}(w)\\
 & =O(h^{q}).
\end{align*}
As long as $q$ can be chosen large enough so that $Nh^{2q}\rightarrow0$
while $Nh^{p}\geq C>0,$ the bias term $T_{4}$ will be asymptotically
negligible and the density estimator $\hat{f}_{W}(w)$ should still
be asymptotically normal with asymptotic distribution of the same
form derived above.

Finally, a particularly useful extension of the kernel estimation
approach for dyadic data would be to estimation of the conditional
expectation of one dyadic variable $Y_{ij}$ conditional on the value
$w$ of another dyadic variable $W_{ij},$ i.e., estimation of 
\[
g(w)\equiv\mathbb{E}[Y_{ij}|W_{ij}=w]
\]
when the vector $W_{ij}$ has $p$ jointly-continuously distributed
components conditional upon any remaining components. Here the Nadaraya-Watson
kernel regression estimator \citep{Nadaraya_TPA64,Watson_Sankhya64}
would be defined as 
\[
\hat{g}(w)\equiv\frac{\sum_{i\neq j}K\left(\frac{w-W_{ij}}{h}\right)Y_{ij}}{\sum_{i\neq j}K\left(\frac{w-W_{ij}}{h}\right)},
\]
and the model for the dependent variable $Y_{ij}$ would be analogous
to that for $W_{ij},$ with 
\begin{align*}
Y_{ij} & =Y(A_{i},B_{j},U_{ij})\\
W_{ij} & =W(A_{i},B_{j},V_{ij})
\end{align*}
in the directed case (and $B_{j}\equiv A_{j}$ for undirected data),
with $\{A_{i}\},$ $\{B_{j}\},$ and $\{(U_{i},V_{ij})\}$ assumed
mutually independent and i.i.d. The large-sample theory would treat
the numerator of $\hat{g}(w)$ similarly to that for the denominator
(which is proportional to the kernel density estimator $\hat{f}_{W}(w)$);
our initial calculations for undirected data with a scalar, continuously-distributed
regressor $W_{ij}$ yield 
\[
\sqrt{N}\left(\hat{g}(w)-g(w)\right)\overset{D}{\rightarrow}\mathcal{N}(0,4\Gamma_{1}(w)),
\]
when $Nh^{p}\rightarrow\infty$ and $Nh^{4}\rightarrow0$, where 
\[
\Gamma_{1}(w)\equiv\mathbb{V}\left(\frac{\mathbb{E}[Y_{ij}|A_{i},W_{ij}=w]\cdot f_{W|A}(w|A_{i})}{f_{W}(w)}\right).
\]
If this calculation is correct, then, like the density estimator $\hat{f}_{W}(w)$
the rate of convergence for the estimator $\hat{g}(w)$ of the conditional
mean $g(w)$ would be the same as the rate for the estimator $\hat{\mu}_{Y}=\bar{Y}$
of the unconditional expectation $\mu_{y}=\mathbb{E}[Y_{ij}]=\mathbb{E}[g(W_{ij})],$
in contrast to the estimation using i.i.d. (monadic) data. We intend
to verify these calculations and derive the other extensions in future
work.

\appendix

\section{\label{app: Proofs}Proofs}

\subsubsection*{Derivation of bias expression, equation \eqref{eq: bias} of the
main text}

Under the conditions imposed in the main text, the expected value
of $\hat{f}_{W}(w)$ is 
\begin{align*}
\mathbb{E}\left[\hat{f}_{W}(w)\right] & =\mathbb{E}\left[\frac{1}{h}K\left(\frac{w-W_{12}}{h}\right)\right]\\
 & =\mathbb{E}\left[\int\frac{1}{h}K\left(\frac{w-s}{h}\right)f_{W}(s)\mathrm{d}s\right]\\
 & =\int K\left(u\right)f_{W}(w-hu)\mathrm{d}u\\
 & =f_{W}(w)+\frac{h^{2}}{2}\frac{\partial^{2}f_{W}(w)}{\partial w{}^{2}}\int u^{2}K\left(u\right)\mathrm{d}u+o(h^{2})\\
 & \equiv f_{W}(w)+h^{2}B(w)+o(h^{2}).
\end{align*}
The first line in this calculation follows from the fact that $W_{ij}$
is identically distributed for all $i,j$, the third line uses the
change-of-variables $s=w-hu,$ and the fourth line follows from a
second-order Taylor's expansion of $f_{W}(w-hu)$ around $h=0$ and
the fact that 
\[
\int u\cdot K(u)\mathrm{d}u=0
\]
because $K(u)=K(-u).$

\subsubsection*{Demonstration of asymptotic negligibility of $T_{2}$ and $T_{4}$}

Equation \eqref{eq: proj_error2}, which defines $T_{2}$, involves
averages of the random variables 
\begin{align*}
\mathbb{E}[K_{ij}|A_{i},A_{j}] & =\int\frac{1}{h}K\left(\frac{w-s}{h}\right)f_{W|AA}(s|A_{i},A_{j})\mathrm{d}s\\
 & =\int K\left(u\right)f_{W|AA}(w-hu|A_{i},A_{j})\mathrm{d}u
\end{align*}
and 
\begin{align*}
\mathbb{E}[K_{ij}|A_{i}] & =\int\frac{1}{h}K\left(\frac{w-s}{h}\right)f_{W|A}(s|A_{i})\mathrm{d}s\\
 & =\int K\left(u\right)f_{W|A}(w-hu|A_{i})\mathrm{d}u
\end{align*}
which are both assumed bounded, so $T_{2}$ can be written, after
some re-arrangement, as the degenerate second-order U-statistic, 
\[
T_{2}=\frac{1}{n}\sum_{i<j}\left(\mathbb{E}[K_{ij}|A_{i},A_{j}]-\mathbb{E}[K_{ij}|A_{i}]-\mathbb{E}[K_{ij}|A_{j}]+\mathbb{E}[K_{ij}])\right)
\]
with all summands uncorrelated. This implies, squaring and taking
expectations, that 
\begin{align*}
\mathbb{E}[T_{2}^{2}] & =\frac{1}{n^{2}}\sum_{i<j}\mathbb{E}[\left(\mathbb{E}[K_{ij}|A_{i},A_{j}]-\mathbb{E}[K_{ij}|A_{i}]-\mathbb{E}[K_{ij}|A_{j}]+\mathbb{E}[K_{ij}]\right)^{2}]\\
 & \leq\frac{5}{n}\mathbb{E}[(\mathbb{E}[K_{ij}|A_{i},A_{j}])^{2}\\
 & =O\left(\frac{1}{n}\right),
\end{align*}
so 
\[
T_{2}=O_{p}\left(\frac{1}{\sqrt{n}}\right)=O_{p}\left(\frac{1}{N}\right).
\]
Turning to the fourth term, defined in equation \eqref{eq: bias_term},
we demonstrated in Section \ref{sec: MSE} that 
\[
T_{4}=h^{2}B(w)+o(h^{2})=O(h^{2}).
\]

\subsubsection*{Demonstration of consistency of $\hat{\Omega}_{2}\left(w\right)$,
equation \eqref{eq: consistency_of_OMEGA2_tilde} of the main text.}

To show result \eqref{eq: consistency_of_OMEGA2_tilde} of the main
text, we start by showing asymptotic unbiasedness of $\tilde{\Omega}_{2}(w)$
for $\Omega_{2}(w)$. The expected value of the summands in \eqref{eq: OMEGA2_tilde}
equal 
\begin{align*}
\mathbb{E}\left[(K_{12})^{2}\right] & =\frac{1}{h}\int[K\left(u\right)]^{2}f_{W}(w-hu)\mathrm{d}u\\
 & =\frac{f_{W}(w)}{h}\cdot\int[K\left(u\right)]^{2}\mathrm{d}u+O(1)\\
 & \equiv\frac{1}{h}\Omega_{2}(w)+O(1)\\
 & =O\left(\frac{1}{h}\right),
\end{align*}
from which asymptotic unbiasedness follows, since: 
\begin{align*}
\mathbb{E}\left[\tilde{\Omega}_{2}(w)\right] & =h\left[\frac{1}{h}\Omega_{2}(w)+O(1)\right]\\
 & =\Omega_{2}(w)+o(1).
\end{align*}
Following the same logic used to calculate the variance of $\hat{f}_{W}(w),$
we calculate the variance of $\tilde{\Omega}_{2}(w)$ as 
\begin{align*}
\mathbb{V}\left(\tilde{\Omega}_{2}(w)\right) & =\mathbb{V}\left(\frac{h}{n}\sum_{i<j}K_{ij}^{2}\right)\\
 & =\left(\frac{h}{n}\right)^{2}\sum_{i<j}\sum_{k<l}\mathbb{C}(K_{ij}^{2},K_{kl}^{2})\\
 & =\frac{h^{2}}{n}\left[\mathbb{V}(K_{12}^{2})+2(N-2)\cdot\mathbb{C}(K_{12}^{2},K_{13}^{2})\right].
\end{align*}
The first term in this expression depends upon 
\begin{align*}
\mathbb{V}(K_{12}^{2}) & =\mathbb{E}\left[K_{12}^{4}\right]-\mathbb{E}\left[K_{12}^{2}\right]^{2}\\
 & =\frac{f_{W}(w)}{h^{3}}\cdot\int[K\left(u\right)]^{4}\mathrm{d}u+O\left(\frac{1}{h^{2}}\right)-\mathbb{E}\left[K_{12}^{2}\right]^{2}\\
 & =O\left(\frac{1}{h^{3}}\right),
\end{align*}
while the second involves 
\begin{align*}
\mathbb{C}(K_{12}^{2},K_{13}^{2}) & =\mathbb{E}[K_{12}^{2}K_{13}^{2}]-\mathbb{E}\left[K_{12}^{2}\right]^{2}\\
 & =\frac{1}{h^{2}}\mathbb{E}\left[\int\left[K\left(u_{1}\right)\right]^{2}f_{W|A}(w-hu_{1}|A_{1})\mathrm{d}u_{1}\right.\\
 & \cdot\left.\int\left[K\left(u_{2}\right)\right]^{2}f_{W|A}(w-hu_{2}|A_{1})\mathrm{d}u_{2}\right]-\mathbb{E}\left[K_{12}^{2}\right]^{2}\\
 & =O\left(\frac{1}{h^{2}}\right).
\end{align*}
Putting things together we have that 
\begin{align*}
\mathbb{V}\left(\tilde{\Omega}(w)\right) & =\frac{h^{2}}{n}\left[\mathbb{V}(K_{12}^{2})+2(N-2)\cdot\mathbb{C}(K_{12}^{2},K_{13}^{2})\right]\\
 & =\frac{h^{2}}{n}\left[O\left(\frac{1}{h^{3}}\right)+2(N-2)\cdot O\left(\frac{1}{h^{2}}\right)\right]\\
 & =O\left(\frac{1}{nh}\right)+O\left(\frac{1}{N}\right)\\
 & =o(1),
\end{align*}
which, with convergence of the bias of $\tilde{\Omega}_{2}(w)$ to
zero, establishes \eqref{eq: consistency_of_OMEGA2_tilde} of the
main text.

\subsubsection*{Demonstration of consistency of $\hat{\Omega}_{1}(w)$, equation
\eqref{eq: consistency_of_OMEGA1_hat} of the main text.}

Since $\hat{f}_{W}(w)$ is consistent if $Nh^{4}\rightarrow0$ and
$Nh\geq C>0,$ consistency of $\hat{\Omega}_{1}(w)$ depends on the
consistency of 
\[
\hat{\mathbb{E}}[K_{12}K_{13}]\equiv\left(\begin{array}{c}
N\\
3
\end{array}\right)^{-1}\sum_{i<j<k}S_{ijk}
\]
for $\lim_{N\rightarrow\infty}\mathbb{E}[K_{12}K_{13}].$ By the fact
that $K_{ij}=K_{ji},$ the expected value of $\hat{\mathbb{E}}[K_{12}K_{13}]$
is 
\begin{align*}
\mathbb{E}[S_{ijk}] & =\mathbb{E}\left[\frac{1}{3}\left(K_{ij}K_{ik}+K_{ij}K_{jk}+K_{ik}K_{jk}\right)\right]\\
 & =\mathbb{E}\left[K_{12}K_{13}\right]\\
 & =\mathbb{E}\left[\int\left[K\left(u_{1}\right)\right]f_{W|A}(w-hu_{1}|A_{1})\mathrm{d}u_{1}\right.\\
 & \cdot\left.\int\left[K\left(u_{2}\right)\right]f_{W|A}(w-hu_{2}|A_{1})\mathrm{d}u_{2}\right]\\
 & =\mathbb{E}\left[f_{W|A}(w|A_{1})^{2}\right]+o(1)
\end{align*}
from the calculations in Section \ref{sec: MSE} above. To bound the
variance of $\hat{\mathbb{E}}[K_{12}K_{13}],$ we note that, although
$\hat{\mathbb{E}}[K_{12}K_{13}]$ is not a U-statistic, it can be
approximated by the third-order U-statistic 
\[
U_{N}\equiv\left(\begin{array}{c}
N\\
3
\end{array}\right)^{-1}\sum_{i<j<k}p_{N}(A_{i},A_{j},A_{k}),
\]
where the kernel $p_{N}(\cdot)$ is 
\begin{align*}
p_{N}(A_{i},A_{j},A_{k}) & =\mathbb{E}[S_{ijk}|A_{i},A_{j},A_{k}]\\
 & =\frac{1}{3}\left(\kappa_{ijk}+\kappa_{jik}+\kappa_{kij}\right),
\end{align*}
for 
\begin{align*}
\kappa_{ijk} & \equiv\mathbb{E}[K_{ij}K_{ik}|A_{i},A_{j},A_{k}]\\
 & =\int\int\frac{1}{h^{2}}\left[K\left(\frac{w-s_{1}}{h}\right)\right]\cdot\left[K\left(\frac{w-s_{2}}{h}\right)\right]\\
 & \cdot f_{W|AA}(s_{1}|A_{i},A_{j})f_{W|AA}(s_{2}|A_{i},A_{k})\mathrm{d}s_{1}\mathrm{d}s_{2}\\
 & =\int\left[K\left(u_{1}\right)\right]f_{W|AA}(w-hu_{1}|A_{i},A_{j})\mathrm{d}u_{1}\\
 & \cdot\int\left[K\left(u_{2}\right)\right]f_{W|AA}(w-hu_{2}|A_{i},A_{k})\mathrm{d}u_{2}.
\end{align*}
The difference between $\hat{\mathbb{E}}[K_{12}K_{13}]$ and $U_{N}$
is 
\[
\hat{\mathbb{E}}[K_{12}K_{13}]-U_{N}\equiv\left(\begin{array}{c}
N\\
3
\end{array}\right)^{-1}\sum_{i<j<k}\left(S_{ijk}-\mathbb{E}\left[S_{ijk}|A_{i},A_{j},A_{k}\right]\right),
\]
and the independence of $\{V_{ij}\}$ and $\{A_{i}\}$ across all
$i$ and $j$ implies that all terms in this summation have expectation
zero and are mutually uncorrelated with common second moment, so that
\begin{align*}
\mathbb{E}\left[\left(\hat{\mathbb{E}}[K_{12}K_{13}]-U_{N}\right)^{2}\right] & \equiv\left(\begin{array}{c}
N\\
3
\end{array}\right)^{-1}\mathbb{E}[\left(S_{123}-\mathbb{E}\left[S_{123}|A_{1},A_{2},A_{3}\right]\right)^{2}]\\
 & \leq\left(\begin{array}{c}
N\\
3
\end{array}\right)^{-1}\mathbb{E}[S_{123}^{2}].
\end{align*}
But 
\begin{align*}
\mathbb{E}[\left(S_{123}\right)^{2}] & =\mathbb{E}\left[\frac{1}{3}\left(K_{12}K_{13}+K_{12}K_{23}+K_{13}K_{23}\right)\right]^{2}\\
 & =\frac{1}{9}\left(3\mathbb{E}\left[\left(K_{12}K_{13}\right)^{2}\right]+6\mathbb{E}[K_{12}^{2}K_{13}K_{23}]\right),
\end{align*}
where 
\[
\mathbb{E}\left[\left(K_{12}K_{13}\right)^{2}\right]=O\left(\frac{1}{h^{2}}\right),
\]
from previous calculations demonstrating consistency of $\Omega_{2}\left(w\right)$,
and 
\begin{align*}
\mathbb{E}[K_{12}^{2}K_{13}K_{23}] & =\mathbb{E}\left[\int\int\int\frac{1}{h^{4}}\left[K\left(\frac{w-s_{1}}{h}\right)\right]^{2}\cdot\left[K\left(\frac{w-s_{2}}{h}\right)\right]\cdot\left[K\left(\frac{w-s_{2}}{h}\right)\right]\right.\\
 & \cdot\left.f_{W|AA}(s_{1}|A_{1},A_{2})f_{W|AA}(s_{2}|A_{1},A_{3})f_{W|AA}(s_{2}|A_{1},A_{3})\mathrm{d}s_{1}\mathrm{d}s_{2}\mathrm{d}s_{3}\right]\\
 & =\frac{1}{h}\mathbb{E}\left[\int\left[K\left(u_{1}\right)\right]^{2}f_{W|AA}(w-hu_{1}|A_{1},A_{2})\mathrm{d}u_{1}\right.\\
 & \cdot\left.\int K\left(u_{2}\right)f_{W|AA}(w-hu_{2}|A_{1},A_{3})\mathrm{d}u_{2}\right]\\
 & \cdot\left.\int K\left(u_{2}\right)f_{W|AA}(w-hu_{2}|A_{1},A_{3})\mathrm{d}u_{2}\right]\\
 & =O\left(\frac{1}{h}\right).
\end{align*}
These results generate the inequality 
\begin{align*}
\mathbb{E}\left[\left(\hat{\mathbb{E}}[K_{12}K_{13}]-U_{N}\right)^{2}\right] & \leq\left(\begin{array}{c}
N\\
3
\end{array}\right)^{-1}\mathbb{E}[\left(S_{123}\right)^{2}]\\
 & =\left(\begin{array}{c}
N\\
3
\end{array}\right)^{-1}\left(O\left(\frac{1}{h^{2}}\right)+O\left(\frac{1}{h}\right)\right)\\
 & =O\left(\frac{1}{N(Nh)^{2}}\right)\\
 & =o(1).
\end{align*}
Finally, we note that $U_{N}$ is a third-order ``smoothed'' U-statistic
with kernel 
\[
p_{N}(A_{i},A_{j},A_{k})=\frac{1}{3}\left(\kappa_{ijk}+\kappa_{jik}+\kappa_{kij}\right)
\]
satisfying 
\[
\mathbb{E}\left[\left(p_{N}(A_{i},A_{j},A_{k})\right)^{2}\right]=O(1)
\]
by the assumed boundedness of $K(u)$ and the conditional density
$f_{W|AA}(w|A_{i},A_{j}).$ Therefore, by Lemma A.3 of \citet{Ahn_Powell_JOE93},
\begin{align*}
U_{n}-\mathbb{E}[U_{N}] & =U_{N}-\mathbb{E}[S_{ijl}]\\
 & =U_{N}-\mathbb{E}\left[f_{W|A}(w|A_{1})\right]^{2}+o(1)\\
 & =o_{p}(1).
\end{align*}
Finally, combining all the previous calculations, we get 
\begin{align*}
\hat{\Omega}_{1}(w) & =\hat{\mathbb{E}}[K_{12}K_{13}]-\left(\hat{f}_{W}(w)\right)^{2}\\
 & =\left(\hat{\mathbb{E}}[K_{12}K_{13}]-U_{N}\right)+(U_{N}-\mathbb{E}\left[f_{W|A}(w|A_{1})\right]^{2})+\mathbb{E}\left[f_{W|A}(w|A_{1})\right]^{2}\\
 & -\left(\left(\hat{f}_{W}(w)\right)^{2}-\left(f_{W}(w)\right)^{2}\right)-\left(f_{W}(w)\right)^{2}\\
 & =E\left[f_{W|A}(w|A_{1})\right]^{2}-\left(f_{W}(w)\right)^{2}+o_{p}(1)\\
 & \equiv\Omega_{1}(w)+o_{p}(1),
\end{align*}
as claimed.

 \bibliographystyle{apalike2}
\bibliography{../Networks_Book/Reference_BibTex/Networks_References}

\end{document}